\documentclass[12pt, oneside,reqno]{amsproc}   	% use "amsart" instead of "article" for AMSLaTeX format
\usepackage{geometry}                		% See geometry.pdf to learn the layout options. There are lots.
\usepackage{graphicx}				% Use pdf, png, jpg, or eps with pdflatex; use eps in DVI mode
\usepackage{amssymb}
\usepackage{amsmath}
\usepackage[all]{xy}
\usepackage{mathtools}
\usepackage{tikz-cd}
\usepackage{enumitem}

\usepackage{tikz-cd}

\newtheorem{thm}{Theorem}[section]
\newtheorem{lem}[thm]{Lemma}
\newtheorem{cor}[thm]{Corollary}
\newtheorem{prop}[thm]{Proposition}

\newtheorem{defn}{Definition}[section]

\numberwithin{equation}{section}

\renewcommand{\t}{\tau}

\def\Pb{\ifmmode{\Bbb P}\else{$\Bbb P$}\fi}
\def\Z{\ifmmode{\Bbb Z}\else{$\Bbb Z$}\fi}
\def\C{\ifmmode{\Bbb C}\else{$\Bbb C$}\fi}
\def\R{\ifmmode{\Bbb R}\else{$\Bbb R$}\fi}
\def\S{\ifmmode{S^2}\else{$S^2$}\fi}

\def\grad{\operatorname{grad}}

\def\S{\cal S}

\newenvironment{pf}{\paragraph{Proof:}}{\hfill$\square$ \newline}

\begin {document}
	
\title{Entropy and generic mean curvature flow in curved ambient spaces}
\begin{abstract}In the article, we generalize some recent results of Colding and Minicozzi on generic singularities of mean curvature flow to curved ambient spaces. To do so, we make use of a weighted monotonicity formula to derive an ``almost monotonicity" for the entropy upon embedding into $\R^\ell$. We are also lead to study the continuity of the entropy functional in certain cases. \end{abstract}

\author {Alexander Mramor}
\address{Department of Mathematics, University of California Irvine, CA 92617}
\email{mramora@uci.edu}

\maketitle

\setcounter{section}{0}

\section{Introduction.}

In this article we show how some of the theory of \cite{CM} by Colding and Minicozzi, on the study of singularities under the mean curvature flow (MCF) through the entropy functional, may be generalized from flows in $\R^3$ to flows in general curved ambient spaces $N^3$. Our main result, showcasing how their theory may be applied in the curved setting, is the following theorem generalizing Theorem 0.10 in \cite{CM}:
\begin{thm} For any closed embedded surface $M$ contained in a 3 manifold $N^3$, there exists a piece-wise MCF $M_t$ starting at $M$ and defined up to time $t_0<\infty$ where the surfaces become singular. Moreover, $M_t$ can be chosen so that if 
\begin{equation}
\liminf\limits_{t \to \t_0} \frac{diam M_t}{\sqrt{t_0 - t}} < \infty
\end{equation}
then $M_t$ becomes extinct in a round point. 
\end{thm}
See the concluding remarks, section 7 below, for a discussion on how this result relates to preexisting results. Recall a mean curvature flow of a smooth submanifold $M \subset N$ is a map $\mathcal{F}: M \times [0, t_0) \to N$ so that $\frac{d\mathcal{F}}{dt} = -H$, where $H$ is the mean curvature vector - see \cite{Mant} for a nice introduction.  We typically denote (like in the above statement) the image of $\mathcal{F}$ at time $t$ by $M_t$ for notational convenience. We will also often use the Brakke flow definition/language of MCF, see Brakke's thesis \cite{B}. Of course for smooth $M$, the mean curvature flow equation $\frac{d\mathcal{F}}{dt} = -H$ (which happens to be a degenerate parabolic PDE) with initial condition $M$ has a solution (which is to say a MCF out of $M$ exists) for at least some short time. Denoting the maximal time of existence by $t_0$, if $t_0 < \infty$ we call this time the $\textit{singular time}$. As $t \to t_0$ necessarily the second fundamental form of $M$ blows up at some point and we call the point(s) $(x_0, t_0)$ in spacetime where it does $\textit{singular points}$. 
$\medskip$

Singularities are practically inevitable in the mean curvature flow (and indeed are the subject of this paper). To study them one successively rescales around where the singularity is occurring and studies the limit flow $\mathcal{B}_t$ which we define first as follows for a MCF $M_t \subset \R^\ell$, following the exposition in \cite{CM}. Let $(x_j, t_j)$ be a sequence of points in spacetime and $c_j$ a sequence of positive numbers with $c_j \to \infty$. A $\textit{limit flow}$ is a Brakke limit (varifold convergence at each time slice) of the sequence of rescaled flows $t \to c_j(M_{c_j^{-2} t + t_j} - x_j)$. Such limits exist by Brakke's compactness theorem, see \cite{B}. When all $x_j$ are equal to a common point $x_0$ and the $t_j$'s are equal to a common time $t_0$, then the limit flow is called a $\textit{tangent flow}$ at $(x_0, t_0)$. 
$\medskip$

In case $M_t$ is a flow is some ambient manifold $N$, we embed a neighborhood of $N$ about the point we are interested in isometrically into $\R^\ell$ for some $\ell$ using Nash's theorem \cite{N} and proceed as above in the definition of tangent flow (although $M_t$ won't be a mean curvature flow in $\R^\ell$, a point we'll come back to shortly). 
$\medskip$

Returning to the statement of the theorem above, by $M_t$ becoming extinct is a round point we mean that all the tangent flows of $M_t$ at a singular point are shrinking round spheres. Equation (1.1) above implies (see \cite{CM}) that there exists some $D > 0$ so that the tangent flows of $M_t$ have diameter $<D$ for all time slices. Furthermore we also assume our flows are smooth up to $\textit{and including}$ the first singular point as Colding and Minicozzi do in that we assume for simplicity every tangent flow is multiplicity one (see also \cite{I}).
$\medskip$

$\textbf{Acknowledgements:}$ The author is extremely grateful to his advisor, Richard Schoen, for his advice and patience. The author is indebted to Jonathan Zhu for informing him on the weighted monotonicity formulas in White's paper \cite{W}. The author also thanks Chao Li and Christos Mantoulidis for making helpful suggestions on clarifying the statement of theorem 1.1.

\section{The first technical hurdle: failure of Huisken monotonicity.} 

Consider for now a mean curvature flow $M_t^n$ in $\R^N$. Denote by $\Phi_{x_0, t_0}$ the backward heat kernel at $(x_0, t_0)$, that is: 
\begin{equation}
\Phi_{x_0, t_0}(x,t) = \frac{1}{(4 \pi (t_0 - t))^{n/2}} \cdot \text{exp}\left ( - \frac{|x|^2}{4 (t_0 - t)} \right ), \text{ }t < t_0.
\end{equation} 
Then Huisken's montonicity (theorem 3.1 in \cite{H}) says the integral of $\Phi_{x_0, t_0}$ is nonincreasing under the flow; more precisely 
\begin{thm}($\textbf{Huisken monotonicity}$) If $M_t$ is a surface flowing by the mean curvature flow for $t < t_0$, then we have the formula
\begin{equation}
\frac{d}{dt} \int_{M_t} \Phi_{x_0, t_0}(x, t) d\mu_t = - \int_{M_t} \Phi_{x_0, t_0}(x,t) \left | H + \frac{1}{2(t_0 - t)} \mathcal{F}^\perp \right |^2 d\mu_t
\end{equation}
\end{thm}
Flows that make this derivative zero (so that $ H + \frac{1}{2(t_0 - t)} \mathcal{F}^\perp = 0$) are called $\textit{self shrinkers}$ and are important below. This is monotonicity is also important in the proof of theorem 0.10 (the Euclidean case) in giving us local volume bounds along the flow depending on initial data and letting us use Brakke's regularity theorem in a critical step (see lemmas 5.1 and 5.2 below). For general ambient manifolds though we don't have such a clean equation for the backwards heat kernel and don't have such a useful quantity right away. 
$\medskip$

The natural thing to do then when the ambient space isn't $\R^N$ is to try to isometrically embed, at least locally, our ambient space into $\R^\ell$ for some $\ell$ by Nash's embedding theorem. When we do this though the flow of $M$ (as a flow in $\R^\ell$) is not a mean curvature flow but instead involves forcing terms from the curvature of the embedding of $N$ in $\R^\ell$. Namely, denoting by $\mathcal{F}: M \times [0, t_0)$ the mean curvature flow of $M$ in $N$ as a flow in $\R^\ell$ we see that $\frac{d\mathcal{F}}{dt} = -H \nu - \text{trace}(B(x) \mid T_xM_t)$ where $B$ is the second fundamental form of $N$ in $\R^\ell$. One finds, just naively calculating the time derivative of $\int \Phi_{x_0, t_0}$, extra terms show up due to the forcing term $P=- \text{trace}(B(x) \mid T_xM_t)$ . Indeed for any function $\psi: \R^\ell \to \R$: 
\begin{equation}
\frac{d}{dt} \int_{M_t} \psi = \int_{M_t} (- \psi |H|^2 + \nabla \psi \cdot H + (\nabla \psi - \psi H \nu) \cdot P^\perp) 
\end{equation}
The idea then is to introduce a new quantity that is monotone under the forced flow that stays close enough to $\int \Phi_{x_0, t_0}$ to ``carry" it along, giving us an almost monotonicity of the quantities we are actually interested about in the necessary circumstances. 

\section{A weighted monotonicity formula for forced flows.} 

By scaling the Huisken weighted volume above by an appropriate weight, we find a monotone quantity under the forced flow (hence the namesake weighted monotonicity formula) that fits our needs. We follow the discussion of it in section 11 of White's stratification paper \cite{W}, where for example the weighted monotonicity formula is used to show stratification results by White for the singular set of the mean curvature flow in $\R^n$ are also valid in curved ambient spaces. To begin we define a $K$-almost Brakke flow: 
\begin{defn} A one-parameter family $M = \{(t, \mu_t) \mid a < t < b \}$ of radon measures in $U \subset \R^\ell$ is a $K\textbf{-almost Brakke flow}$ provided
\begin{enumerate}
\item For almost every $t$, $\mu_t$ is the radon measure associated with an integer multiplicity rectifiable varifold.
\smallskip

\item For every nonnegative compactly supported $C^1$ function $\phi$ on $U$,
\begin{equation} 
\overline{D}_t \int \phi d\mu_t \leq \int  (- \phi |H|^2 + \nabla \phi \cdot H + K (\nabla \phi - \phi H \nu) )  d\mu_t
\end{equation}
where $\overline{D}f(t) := \lim\sup\limits_{h \to 0} \frac{f(t + h) - f(t)}{h}$
\end{enumerate}
\end{defn}
Note, that for any smooth (local) isometric embedding of $U \subset N^3$ into $\R^\ell$,  the forcing term $P$ will be bounded by some $K$ depending on the second fundamental of the embedding so the mean curvature flow of a hypersurface $M \subset U$ in $N$ will be a $K$-almost flow in $\R^\ell$. With this terminology we give the anticipated weighted monotonicity formula, which in following with \cite{W} we first present in its most general (i.e localized) form: 
\begin{lem} ($\textbf{weighted monotonicity - general form}$) Let $\mathcal{M} = \{(t, \mu_t) \mid a < t < b\}$ be an $m$-dimensional $K$-almost Brakke flow in $U \subset \R^\ell$ with $\mu_t(U) \leq \Lambda < \infty$ for all $t$. Suppose $s \in (a,b)$ and $\mathbf{B}(y, 2r) \subset U$, and let $\psi: \mathbf{B}(2y,r) \to [0,1]$ be a $C^2$ function that is 1 in $\mathbf{B}(y,r)$ and satisfies the bound:
\begin{equation}
r| D \psi| + r^2|D^2 \psi| \leq b
\end{equation}
Then the function 
\begin{equation} 
J_{y,s}(t) = e^{K^2(s-t)/2}u_{y,s} + \left ( \frac{c_m(1 + b) \Lambda}r^{m+2} \right) \frac{e^{K^2(s-t)/2} - 1}{K^2/2}
\end{equation}
is non-increasing on the interval $\max \{ s - r^2, a \} \leq t < s$ where 
\begin{equation}
\rho = \Phi_{y, s}(x,t) \text{, } u_{y,s} = \int \psi \rho
\end{equation}
Indeed,
\begin{equation} 
J(t_2) - J(t_1) \leq - \frac{1}{2} e^{K^2(s - t_2)/2} \int_{t_1}^{t_2} \int \psi \rho \left | J + \frac{(x-y)^\perp}{2(s-t)} - \frac{(D \psi)^\perp}{\psi} \right |^2 d \mu_t  dt
\end{equation}
\end{lem}
Of course since mass decreases under the flow and our flows are compact smooth surfaces, we automatically get mass bounds $\Lambda < \infty$ for any open set $U \subset \R^\ell$. Note that if we let $\psi = 1$ on all of $\R^\ell$, we may take bound $b$ above to be zero and $r$ as large as we want, deriving all we really need (we adapt the definition of $u$ below): 
\begin{cor} ($\textbf{weighted monotonicity - simpler form}$) Let $\mathcal{M} = \{(t, \mu_t) \mid a < t < b\}$ be an $m$-dimensional $K$-almost Brakke flow in $\R^\ell$ with uniformly bounded mass. 
Then the function 
\begin{equation} 
J(t) = e^{K^2(s-t)/2} \int \Phi_{y,s}(x,t)  = e^{K^2(s-t)/2} u_{y,s}(x,s)
\end{equation}
is non-increasing on the interval $a \leq t < s$. Indeed for $t_1 < t_2 < s$:
\begin{equation} 
J(t_2) - J(t_1) \leq - \frac{1}{2} e^{K^2(s - t_2)/2} \int_{t_1}^{t_2} \int \psi \rho \left | J + \frac{(x-y)^\perp}{2(s-t)}  \right |^2 d \mu_t x dt
\end{equation}
\end{cor}

Even though we want to use the second version of the result above we might as well keep the notation $u_{y,s} =  \int \Phi_{y,s}(x,t)$, so that $J =e^{K^2(s-t)/2} u_{y,s}$. Before moving to the ``almost monotonicity" statement first we record an important corollary of lemma 3.1 (besides White's stratification results); if $M_t \subset \R^\ell$ is a $K$-almost Brakke flow, then dilation of $\R^\ell$ by $S$ is a $K/S$-almost Brakke flow, implying by Brakke compactness (which is also true for $K$-almost Brakke flows) that the tangent flow is a regular Brakke flow in $\R^\ell$. Even more, the weighted montonicity implies that Huisken's density is upper-semicontinuous so arguing as in \cite{I} the tangent flows are actually self-shrinkers in $\R^\ell$:
\begin{cor} Tangnent flows to $K$-almost Brakke flows are ordinary Brakke flows, furthermore they are self shrinkers. \end{cor}
$\medskip$

In fact, we immediately see something a bit more is true that we'll want for the sequel. It's an important observation because we will want to import the theory for hypersurface self-shrinkers in $\R^3$ from \cite{CM}; it is not good enough to merely know the tangent flows are surfaces in $\R^\ell$.  We see it's true of course since $N$ is 3 dimensional and rescaling $\R^\ell$ in the tangent flow ``flattens out" $N$: 
\begin{cor} Suppose that $M \times [0,t_0)$ is a flow contained in an open region of $N^3$ that can be embedded into $\R^\ell$ for some $\ell$. so is a $K$-almost Brakke flow for some $K$. Then a tangent flow to $M$ at $x \in N^3$ is an ordinary Brakke flow in $T_x N \cong \R^3 \subset \R^\ell$, where the inclusion $\R^3 \subset \R^\ell$ is flat. 
\end{cor}

Now note that $u_{y,s}$ is indeed Huisken's weighted volume from before. We are interested then when $J$ is close to $u$ and how close it is; to start, note that the monotonicity of $J$ implies that $u$ is uniformly bounded on any finite time interval, even though $u$ itself might not be monotone:
\begin{lem} Suppose that $M_t$ is a $K$-almost Brakke flow and that $u(x,t_0)$ is bounded by $C_1$ for all $x$ at $t_0$. Then if $t - t_0 < \tau < \infty$ there is a constant $C$ so that, $u(x,t) < C$. \end{lem}
\begin{pf} Note that there is a constant $\sigma > 0$ so that for all $t > t_0$, $t - t_0 < \tau$, $e^{K^2(s -t)} > \sigma > 0$. Also, since $J$ is monotone decreasing, $J(t_0) > J(t)$ for any $t$ in this time interval. Since $J = e^{K^2(s-t)}u$ then, $u(t) < J(t_0)/ \sigma < e^{K^2(s - t_0)}C/\sigma$, so we may take $C = e^{K^2(s - t_0)}C_1/\sigma$. \end{pf}
$\medskip$

Note that $|u - J| =|u(1 - e^{K^2(s-t)})| \leq C|1 - e^{K^2(s-t)}|$, using the bound from above Taylor expanding we see that $1 - e^{K^2(s-t)} \leq K^2(s-t)e^{K^2(s-t)}$. The monotonicity of $J$ then gives us:
\begin{prop}($\textbf{almost monotonicity of }\mathbf{u_{y,s}}$) Suppose that $M \times [0,T)$ is a $K$-almost Brakke flow in $\R^\ell$. Then given a point $(y,s)$ in spacetime, $0 <\tau < \infty$, and $C > 1$,  there exists $\delta > 0$ so that if $K < \delta$, $u_{y,s}(t_2) < u_{y,s}(t_1) + CK^2(t_2 - t_1)$ where $t_1 < t_2 < s$, $t_2 - t_1 < \tau$. 
\end{prop}
Again if we rescale $\R^\ell$ by $S$ then $K$ scales by $K \to K/S$, so we can certainly rescale $\R^\ell$ to make the assumptions of the above applicable, but this also would dilate $\tau$ so isn't something we can immediately do profitably. Below we'll still be able to make use of it when considering tangent flows satisfying (1.1) since it implies upper bounds on singular time (after rescaling). 
$\medskip$

\section{Essential facts from Colding-Minicozzi theory and almost monotonicity of entropy.} 
Our proof of theorem 1.1 rests heavily on concepts and terminology of Colding and Minicozzi so we develop the necessary machinery in this section. Since we embed into $\R^\ell$, we will also need some generalizations of their theory to higher codimension - see \cite{YL} and \cite{And}. Consider a surface $\Sigma^k \subset \R^{\ell}$; then given $x_0 \in \R^{\ell}$ and $t_0 > 0$ define the functional $F_{x_0, t_0}$ by 
\begin{equation}
F_{x_0, t_0}(\Sigma) = \frac{1}{(4 \pi t_0)^{k/2}} \int_\Sigma e^\frac{-|x - x_0|^2}{4t_0} d\mu
\end{equation} 
We see that $F$ is closely related to $u$; indeed
\begin{equation}
F_{x_0, t_0}(\Sigma) = \frac{1}{(4 \pi t_0)^{k/2}} \int_{\Sigma} e^{-\frac{|x - x_0|^2}{4t_0}}  d\mu= \int_{\Sigma} \Phi_{x_0, t_0}(x, 0) d\mu = u_{x_0, t_0}(0)
\end{equation}
($u$ implicitly depends on $\Sigma$.). Note that the functionals $F_{x_0, t_0}$ are essentially weighted volumes. Naturally then we next record the first variation of $F_{x_0,t_0}$ along the flow; the one calculated in \cite{CM} is for codimension 1 flows but since we are considering $K$-almost flows in $\R^\ell$ of some arbitrary codimension we need a more general one that can be found in \cite{YL} (and also \cite{And}). It is of course important in the sequel:
\begin{lem} (Theorem 1 in \cite{YL}) Let $\Sigma_s \subset \R^{\ell}$ be an $n$-dimensional complete manifold without boundary which has polynomial volume growth. Suppose that $\Sigma_s \subset \R^\ell$ is a normal variation of $\Sigma$, $x_s, t_s$ are variations of $x_0$ and $t_0$ and $\Sigma_0' = V$, $x_0' = y$ and $t_0' = h$, then $\frac{d}{ds}(F_{x_s, t_s}(\Sigma_s))$ is
\begin{equation} 
\frac{1}{(4 \pi t_0)^{k/2}} \int_\Sigma \left( - \langle V, H + \frac{ x - x_0}{2t_0}  \rangle+ h( \frac{|x- x_0|^2}{4t_0^2} - \frac{n}{2t_0} ) + \frac{\langle x - x_0, y \rangle}{2t_0} \right)  e^{\frac{-|x - x_0|^2}{4t_0}} d \mu
\end{equation}
\end{lem} 
Now we define the entropy, of central importance in Colding-Minicozzi theory:
\begin{defn} The entropy $\lambda$ of $\Sigma$ is the supremum over all $F_{x_0, t_0}$ functionals:
\begin{equation} 
\lambda(\Sigma) = \sup\limits_{x_0, t_0} F_{x_0, t_0}(\Sigma) 
\end{equation} 
\end{defn} 
(Since in the definition of $F$ functional $t_0 > 0$, above points range over $(x_0, t_0)$ where $t_0 > 0$.) Note that $\lambda$ is non-negative and invariant under dilations, rotations, or translations of $\Sigma$. Also, for the mean curvature flow $\Sigma_t$ of $\Sigma$, $\lambda(\Sigma_t)$ is nonincreasing under the flow. Of course though you can define $F$ and hence $\lambda$ along a flow of hypersurfaces in $\R^\ell$ not flowing by MCF, for example forced mean curvature flows which is the point of introducing $J$. Before moving further let's list some slightly strengthened properties on $F$ and the relationship of $\lambda$ and $F$ one finds for the hypersurface case in \cite{CM} we'll need later on (we generalized the codimension of $\Sigma$). It is important to note that monotonicity isn't needed in this lemma:
\begin{lem}(generalized lemma 7.2 in \cite{CM}) If $\Sigma^k \subset \R^{\ell}$ is a smooth complete embedded hypersurface without boundary and with polynomial volume growth, then
\begin{enumerate}
\item $F_{x_0, t_0}(\Sigma)$ is a smooth function of $x_0$ and $t_0$ on $\R^{n+1} \times (0, \infty)$.
\item Given any $t_0 > 0$ and any $x_0$, we have $\partial_{t_0} F_{x_0, t_0}(\Sigma) \geq - \frac{\lambda(\Sigma)}{4} \sup_\Sigma |H|^2$
\item For each $x_0$, $\lim_{t_0 \to 0} F_{x_0, t_0}(\Sigma)$ is 1 if $x_0 \in \Sigma$ and is $0$ otherwise. 
\item If $\Sigma$ is closed, then $\lambda(\Sigma) < \infty$
\end{enumerate}
\end{lem}
\begin{pf} Statements (1), (3), and (4) are clear but (2) requires some more work. Without loss of generality $x_0 = 0$. The first variation formula above implies:
\begin{equation}
\partial_{t_0}F_{0, t_0}(\Sigma) = \frac{1}{(4\pi t_0)^{n/2}} \int_\Sigma \frac{|x|^2 - 2nt_0}{4t_0^2} e^{\frac{-|x|^2}{4t_0}}
\end{equation}
Since $\Delta |x|^2 = 2n - \langle x, H \rangle$ and $\Delta e^f = e^f(\Delta f + | \nabla f|^2 )$, we have 
\begin{equation}
\begin{split}
e^{\frac{|x|^2}{4t_0}} \Delta e^{- \frac{|x|^2}{4t_0}} = \frac{ |x^T|^2}{4t_0^2} - \frac{2n}{4t_0} + \frac{\langle x, H \rangle}{2t_0} =\frac{|x|^2 - 2nt_0}{4t_0^2} - \frac{|x^\perp|}{4t_0^2} + |H| \frac{\langle x, \frac{H}{|H|} \rangle}{2t_0} \\  \leq \frac{|x|^2 - 2nt_0}{4t_0^2} + \frac{|H|^2}{4}
\end{split}
\end{equation}
where the inequality used $2ab \leq a^2 + b^2$. Just like in \cite{CM}, since $\Sigma$ has polynomial volume growth and the vector field $\nabla e^{-\frac{|x|^2}{4t_0}}$ decays exponentially, Stokes' theorem gives 
\begin{equation}
\partial_{t_0} F_{0, t_0}(\Sigma) \geq -  \frac{1}{(4\pi t_0)^{n/2}} \int_\Sigma \frac{|H|^2}{4}  e^{-\frac{|x|^2}{4t_0}} \geq -\frac{1}{4} F_{x_0, t_0}(\Sigma) \sup\limits_{\Sigma} |H|^2 \geq - \frac{\lambda(\Sigma)}{4} \sup\limits_{\Sigma} |H|^2
\end{equation}
Showing (2)
\end{pf}
In some cases (relevant to ours in fact), the entropy is even achieved by some $F_{x_0, t_0}$. This is again a lemma whose proof generalizes immediately to higher codimension case; we put off showing a strengthening of it to the next section:
\begin{lem} (lemma 7.7 in \cite{CM}) If $\Sigma \subset \R^{n+1}$ is a smooth closed embedded hypersurface and $\lambda(\Sigma) > 1$, then there exists $x_0 \in \R^{n+1}$ and $t_0 > 0$ so that $\lambda = F_{x_0, t_0}(\Sigma)$. 
\end{lem}
$\medskip$

Why might entropy be useful? For MCF in $\R^\ell$, Huisken monotonicity can be interpreted as saying that for $s < t < t_0$, $F_{x_0, t_0}(M_t) \leq F_{x_0, t_0 + (t - s)}(M_s)$. Indeed, temporarily making explicit $u$ is a function of submanifold: 
\begin{equation}
\begin{split}
F_{x_0, t_0 }(M_t) = u_{x_0, t_0 + t}(M_t)(t) \leq  u_{x_0, t_0 + t}(M_s)(s) \\ =  u_{x_0, t_0 + (t-s)}(M_s)(0) = F_{x_0, t_0 + (t-s)}(M_s)
\end{split}
\end{equation}
where the inequality used Huisken montonicity. Since entropy is defined as a supremum of $F$ functionals over $\{x_0, t_0\}$ where $t_0 > 0$ then we see that for MCF in $\R^\ell$ entropy is monotone decreasing under the flow. So, as an example, if we can understand the singularities of low entropy (or stable entropy), imposing entropy conditions on $M$ would by monotonicity imply what kind of singularities it can have. 
$\medskip$

Of course for $K$-almost Brakke flows we don't have monotonicity, but from this chain of inequalities above though we see we have the following almost monotonicity statement for $F$ using proposition 3.6 above, the almost monotonicity for $u$:
\begin{lem}($\textbf{almost monotonicity of F}$) Suppose that $M \times [0,T)$ is a $K$-almost Brakke flow in $\R^\ell$. Then given a point $(x_0,t_0)$ in spacetime, $0 < \tau < \infty$, and $C > 1$,  there exists $\delta > 0$ so that if $K < \delta$, $F_{x_0, t_0}(M_t) \leq F_{x_0, t_0 + (t - s)}(M_s) + CK^2(t-s)$ where $s < t < t_0$, $t - s < \tau$. 
\end{lem}
This immediately implies that if $M_t$ is an almost Brakke flow $\overline{D}_t \lambda(M_t) < CK^2$ ($\overline{D}_t$ as given above).  Since $\tau < \infty$, after possibly taking $\delta$ smaller (smallness of $K$) we immediately get the following almost monotonicity for entropy:
\begin{prop}($\textbf{almost monotonicity of entropy}$) Suppose that $M \times [0,T)$ is a $K$-almost Brakke flow in $\R^\ell$. Then given $0 < \tau, \epsilon_0 < \infty $, there exists $\delta > 0$ so that if $K < \delta$, $\lambda(M_s) < \lambda(M_t) + \epsilon_0$  for $t < s < T$, $s - t < \tau$. 
\end{prop} 

To use this we will need to eventually answer the question: how small should $\epsilon$ be? The following two statements we'll see later dictate this. Given a constant $D > 0$, let $S_D = S_{g, \overline{\lambda}, D}$ denote the space of all smooth closed embedded self shrinkers in $\R^3$ with genus at most $g$, entropy at most $\overline{\lambda}$, and diameter at most $D$. From \cite{CM}, \cite{CM1} we know:
\begin{prop} (corollary 8.2 of \cite{CM} or \cite{CM1}) For each fixed $D$, the space $S_D$ is compact. Namely, any sequence in $S_D$ has a subsequence that converges uniformly in the $C^k$ topology (any fixed $k$) to a surface in $S_D$.
\end{prop}
An important corollary of which for us is the following
\begin{cor} (corollary 8.4 of \cite{CM}) Given $D > 0$ there exists $\epsilon > 0$ so that if $\Sigma \in S_D$ is not the round sphere, then there is a graph $\Gamma$ over $\Sigma$ with $\lambda(\Gamma) < \lambda(\Sigma) - \epsilon$. 
\end{cor}
Again we emphasize that the tangent flows to $M$ as a $K$-almost Brakke flow in $\R^\ell$ will be ordinary Brakke flows that will lay in $\R^3 \subset \R^\ell$ and hence self shrinkers considered as flows just in $\R^3$. Furthermore from the bound (1.1) above these corollaries will hold without much trouble below. 
$\medskip$

It is also worth pointing out that getting important integral curvature bounds via the genus depends on dimension specific techniques (one could use Gauss-Bonnet, for example) and thus the proofs of proposition 4.6 and hence corollary 4.7 don't carry over in higher dimensions, hence the dimension restriction on $N^3$. 
$\medskip$

\section{Lipschitz continuity of entropy in certain cases.}
Below in the proof of theorem 1.1 we will want to understand how the entropy behaves on a one parameter family $\Sigma_s$, $s \leq 0 \leq 1$ (not necessarily moving by MCF). As Colding and Minicozzi point out in \cite{CM}, the entropy $\lambda(\Sigma_s)$ does not necessarily depend smoothly on $s$; we are only interested though in when we can say it is continuous. We start by proving what one can interpret as a very weak version of Bernstein and Wang's results \cite{BW} (see also \cite{Z} for the extension to higher dimensions), where as a consequence of their work it was shown the entropy of a closed hypersurface is bounded below by the entropy of the round sphere: 
\begin{lem} Let $\Pi_{C,D}$ be the family of compact closed $k$-submanifolds $\Sigma^k$ bounded locally graphically in $C^{3}$ by $C$ and diam$(\Sigma) \leq D$. Then there exists $\sigma > 0$ so that  $\lambda(\Sigma) > 1 + \sigma > 1$ for $\Sigma \in \Pi_{C,D}$.
\end{lem}
\begin{pf} Suppose not. Then there exists a sequence $\Sigma_i \in \Pi_{C,D}$ so that $\lambda(\Sigma_i) < 1+ \frac{1}{i}$. Taking the limit by Arzela-Ascoli, we get a $C^{2,\alpha}$ converging subsequence for some $0 < \alpha < 1$, which we relabel $\Sigma_i$, converging to say $\Sigma$. Since each of the $F_{x_0, t_0}$ is continuous as a function on submanifolds, we see that for each $(x_0, t_0)$, $F_{x_0, t_0}(\Sigma) \leq 1$. Hence the entropy is equal to 1. 
$\medskip$

Now note under mean curvature flow $\Sigma_t$ of $\Sigma$ (as a submanifold of $\R^\ell$) for all $(x_0, t_0)$, $F_{x_0, t_0}(\Sigma_t)$ stays bounded by 1 by Huisken monotonicity so we get curvature bounds on $\Sigma_t$ for all time by Brakke regularity theorem (as in the proof of lemma 6.2 below) so no singularity develops\footnote{This would be good enough for hypersurfaces since there are no compact closed minimal hypersurfaces and so every compact hypersurface must develop a singularity.}. By the curvature bounds on the flow $\Sigma_t$ we may take a subsequential limit along a sequence of times $t_i \to \infty$ to get a limit surface $S$ which we see must be a self shrinker with entropy 1, hence a plane From the rigidity statement for Gaussian density, see proposition 2.10 in \cite{W1} (of course all self shrinkers are ancient flows).
$\medskip$

 But $\Sigma$ has finite volume and of course it remains bounded under the flow so there is no way then plane can arise as a subsequential limit, so we get a contradiction. 

\end{pf}
$\medskip$

We may also attain via a compactness argument:
\begin{lem} Given $\epsilon > 0$, there exists $T > 0$ so that if $\Sigma^k \in \Pi_{C,D}$ and $t_0 < T$ then $F_{x_0, t_0} < 1 + \epsilon$ for any $x_0$. \end{lem}
\begin{pf}
Suppose not. Then for some $\epsilon > 0$ there is a sequence $\Sigma_i \in \Pi_{C,D}$, with corresponding points $(x_i, t_i)$, $t_i \to 0$, so that $F_{x_i, t_i}(\Sigma_i) > 1 + \epsilon$. By Arzela-Ascoli, passing the limit to a $C^2$ graphically converging subsequence, which we relabel back to $\Sigma_i$, with limit say $\Sigma$. We see each of the $F_{x_i, t_i}$ is continuous as a function on submanifolds, and since each of the $\Sigma_i$ have diameter bounded by $D$ so there is a converging subsequence $x_i$ converging to say $x$. Then we see that $\lim\limits_{t_0 \to 0} F_{x, t_0}(\Sigma) > 1 + \epsilon$, which is a contradiction. 
\end{pf}

With this in hand, we can prove the following strengthening of lemma 7.7 in \cite{CM}, as recorded in lemma 4.3 above:
\begin{prop} For $\Sigma \in \Pi_{C,D}$ as defined above, there is a compact set $A \subset \R^\ell \times (0, \infty)$ depending on $C,D$ so that the entropy for $\Sigma \in \Pi_{C,D}$ is achieved in $A$. \end{prop}
\begin{pf} For each fixed $t_0$, it is easy to see that $\lim\limits_{|x_0| \to \infty} F_{x_0, t_0}(\Sigma) = 0$ by the exponential decay of the weight function together with the compactness of $\Sigma$. In particular, for each fixed $t_0 > 0$, the maximum of $F_{x_0, t_0}(\Sigma)$ is achieved at some $x_0$. Moreover, the first variation formula shows (like in codimension 1 case) that this maximum occurs when the weighted integral of $(x - x_0)$ vanishes, but this can only occur when $x_0$ lies in the convex hull of $\Sigma$. It remains to take the supremum of these maxima as we vary $t_0$ - our task then is to show there are constants $0<  T_1 < T_2 < \infty$ so that if $t_0 \not\in [T_1, T_2]$ then for $\Sigma \in \Pi_{C,D}$ one has $F_{x_0, t_0}(\Sigma) < 1 + \sigma$. 
$\medskip$

Using $\sigma$ from lemma 5.1 above , for $\epsilon = \sigma/2$ in lemma 5.2 there is a $T$ so that for all $\Sigma \in \Pi_{C,D}$, when $t_0 > T$ we have $F_{x_0, t_0}(\Sigma) < 1 + \epsilon < 1 + \sigma < \lambda(\Sigma)$, so we take $T_1 = T$. 
$\medskip$

To get $T_2$ first note that for $\Sigma \in \Pi_{C,D}$ that Vol$(\Sigma)$ is universally bounded by say $C_1$ and that 
\begin{equation}
F_{x_0, t_0}(\Sigma) \leq \frac{1}{(4 \pi t_0)^{n/2}} \text{Vol}(\Sigma)
\end{equation} 
hence we easily attain a $T_2$ so that if $t_0 > T_2$, $F_{x_0, t_0}(\Sigma) < 1 + \sigma$.
\end{pf}
$\medskip$

Consider then a smooth one parameter family $\Sigma_s \subset \Pi_{C,D}$ for some $C,D$ and the corresponding manifold $A$; with this in mind we think of the family of $F$ functionals $F_{x_0, t_0}(\Sigma_s)$ as a single function $\mathbf{F}$ on $A \times [0,1]$, associating to each $s \in [0,1] \to \Sigma_s$. The above proposition can be interpreted as saying for every fixed $s$, $\sup \mathbf{F}(\Sigma_s)$ is attained. 
$\medskip$

Since the curvature of the one paramater family $\Sigma_s$ will be bounded along $[0,1]$, the first variation formula lemma 4.1. gives a gradient bound on $\mathbf{F}$ with say $|\nabla F| < C_2$. Slightly modifying the proof of ``Hamilton's trick" to make use of the uniform gradient bound, see lemma 2.1.3 in \cite{Mant}, yields the following: 
\begin{lem}(lemma 2.1.3 in \cite{Mant}) Let $u: M \times [0,1] \to \R$ be a $C^1$ function with $|\grad u| < C$ such that for every time t, there exists a value $\delta > 0$ and a compact subset $K \subset M-\partial M$ such that at every time $t' \in (t- \delta, t + \delta) \cap [0,1]$ the maximum $u_{max}(t') = \max_{p \in M}u(p, t')$ is attained at least at one point of $K$. Then $u_{max}$ is a Lipschitz function in $[0,1]$ with Lipschitz constant $C$. 
\end{lem}
Our manifold $M$ in the above is $A$; by slightly enlarging $A$ above we may ensure that the entropy is attained away from $\partial A$. Hence we derive the following: 
\begin{prop} Suppose that $\Sigma_s^k \subset \R^\ell$, $0 \leq s \leq 1$ is a one-parameter family of closed compact submanifolds bounded locally graphically in $C^{2,\alpha}$. Then $\lambda(\Sigma_s)$ is continuous in $s$.
\end{prop}

\section{Proof of theorem 1.1.}

Consider a hypersurface $M$ flowing to a point in $N^3$ as above. Then there is a time $T > 0$ so that for $t > T$, $M_t$ is contained in a ball of say radius 1 about some point $x_0 \in N^3$. We may isometrically embed $B_{x_0}(1)$ into $\R^\ell$ for some $\ell$, by Nash's embedding theorem; without loss of generality then $N$ is a submanifold of $\R^\ell$ with bounded second fundamental form so that $M$ is a $K$-almost Brakke flow. We start with a couple lemmas, slight modifications of those in \cite{CM}; Huisken monotonicity plays an important role in the original proofs of both of them and so we must modify them to use almost monotonicity. Also, for the second lemma, we have to use the higher codimension first variation formula for $F$ we recorded above. 
$\medskip$

\begin{lem} (modification of lemma 2.9 in \cite{CM}) Let $M_t$ be an n dimensional smooth $K$-almost Brakke flow and choose some $S > 0$. Furthermore suppose $K$ is so small to make the assumptions of proposition 3.6 hold for $ \tau =S, C = 2$. Then given $T>0$, there exists a constant $V = V(\text{Vol}(M_0),S, T) > 0$ so that for all $r < \sqrt{S}$, all $x_0 \in \R^\ell$, and all $T < t$
\begin{equation}
\text{Vol}(B_r(x_0) \cap M_t) \leq (V + 2S)r^n
\end{equation}
\end{lem}
\begin{pf} 
Possibly taking $\delta$ even smaller without loss of generality $K < 1$. For any $t_0 > t$ with $t_0 - t < S$ to be chosen later: 
\begin{equation}
\begin{split}
\frac{1}{(4\pi(t_0 -t))^{n/2}} e^\frac{-1}{4} Vol(B_{\sqrt{t_0 - t}} (x_0) \cap M_t) \leq \frac{1}{(4\pi(t_0 -t))^{n/2}}  \int_{B_{\sqrt{t_0 - t}} (x_0) \cap M_t} e^{\frac{|x - x_0|^2}{4(t - t_0)}} \\ \leq \int_{M_t} \rho_{x_0, t_0}( \cdot, t) = u_{x_0, t_0}(t) \leq u_{x_0, t_0}(0) + CK^2(t_0 -t) \leq \frac{1}{(4\pi T)^{n/2}} \text{Vol}(M_0) + 2S
\end{split}
\end{equation}
Setting $t_0 = t + r^2$ (by assumption, $r < \sqrt{S}$) and multiplying through we get the inequality. 
\end{pf}

\begin{lem} (modification of lemma 8.7 in \cite{CM}) Suppose that $M_t \subset N$,  is a MCF of smooth closed surfaces for $t < 0$ in a smooth manifold $N$ considered embedded in $\R^\ell$ with bounded second fundamental form (so $M_t$ is a $K$-almost Brakke flow in $\R^\ell$ for $t < 0$ for some $K$ depending on the bound). Also suppose $\Sigma_0$ is a closed smooth self shrinker equal to the $t = -1$ time-slice of a multiplicity one tangent flow to $M_t$ at $(0,0) \subset \R^\ell \times \R$. Then we can choose a sequence $s_j > 0$ with $s_j \to 0$ so that 
\begin{equation}
\frac{1}{\sqrt{s_j}} M_{-s_j} \text{ converges in } C^2 \text{ to } \Sigma_0
\end{equation} 
\end{lem}
$\textit{remark:}$ By this convergence in $C^2$ we mean $\frac{1}{\sqrt{s_j}} M_{-s_j}  = \Sigma_0 + X_j$, where $X_j$ is a sequence of vector fields in $N\Sigma_0$ that converges to zero in $C^2$. Also note if $M$ converges to a point we can always arrange it happens at $(0,0)$ by translating. 
$\medskip$

\begin{pf} Fix $\epsilon > 0$ small (to be given by Brakke regularity theorem in \cite{W1}). Since $\Sigma_0$ is a smooth closed embedded surface from lemma 4.2 (1) and (3) above there is $\overline{r} > 0$ so that 
\begin{equation}
\sup\limits_{t_0 \leq \overline{r}} \left( \sup\limits_{x_0 \in \R^3} F_{x_0, t_0}(\Sigma_0) \right) < 1 + \epsilon
\end{equation}
The definition of tangent flows gives a sequence $s_j > 0$ with $s_j \to 0$ so that the rescaled flows $M_t^j = \frac{1}{\sqrt{s_j}} M_{\frac{t}{s_j}}$ converge to the multiplicity-one flow $\sqrt{-t} \Sigma_0$. Let $M_{-1}^j = \frac{1}{\sqrt{s_j}} M_{\frac{-1}{s_j}}$ be the $t = -1$ slice of the $j$-th rescaled flow. We can assume that the $M_{-1}^j$'s converge to $\Sigma_0$ as Radon measures with respect to Hausdorff distance. 
We will use the convergence together with the above to get uniform bounds for the $F$ functionals on the $M_{-1}^j$'s, To do this, define a sequence of functions $g_j$ by 
\begin{equation}
g_j(x_0, t_0) = F_{x_0, t_0}(M_{-1}^j).
\end{equation}
 We will only consider the $g_j$'s on the region $\overline{B} \times [\overline{r}/3, \overline{r}]$ where $B \subset \R^\ell$ is a fixed ball of say radius $2D $ that contains $\Sigma_0$ and all of the $M_{-1}^j$'s. If we have uniform local area bounds for the $M_{-1}^j$ it follows from the first variation formula for $F_{x_0, t_0}$ (lemma 4.1 above) that the $g_j$'s are uniformly Lipschitz in this region with 
 \begin{equation}
 \sup\limits_{\overline{B} \times [\overline{r}/3, \overline{r}]} | \nabla_{x_0, t_0} g_j | \leq C
 \end{equation}
 where $C$ depends on $\overline{r}$, the radius of the ball $B$, and the local area bounds. Of course, since $s_j \to 0$ and $K$ scales by $K \to \sqrt{s_j}K$ there is a $j_0$ so that for $j > j_0$ almost monotonicity holds (let $\tau = \sqrt{2D} + 1$, let $\epsilon_0 = \epsilon$) and hence lemma 5.1 holds so we get the gradient bound in the parabolic cylinder. Since the $M_{-1}^j$'s converge to $\Sigma_0$ as Radon measures and $\Sigma_0$ satisfies (5.4), it follows that 
 \begin{equation} 
 \lim\limits_{j \to \infty} g_j(x_0, t_0) < 1 + \epsilon \text{ for each fixed } (x_0, t_0) \in B \times [\overline{r}/3, \overline{r} ]
 \end{equation}
 Combining this with the derivative estimate and the compactness of $\overline{B} \times [\overline{r}/3, \overline{r} ]$, there exists $j_1 > j_0$ sufficiently large so that for all $j > j_1$ we have 
 \begin{equation}
\sup\limits_{\overline{B} \times [\overline{r}/3, \overline{r}]} F_{x_0, t_0}(M_{-1}^j) = \sup\limits_{\overline{B} \times [\overline{r}/3, \overline{r}]} g_j(x_0, t_0) < 1 + 2 \epsilon
\end{equation}
We see by lemma 4.4, almost-monotonicity of $F$, that after possibly taking $j_0$ larger to ensure in lemma 4.4 that $CK^2 < \epsilon$, for every $-1 < t < 0$ and $j > j_0$:
\begin{equation}
F_{x_0, t_0}(M_t^j) \leq F_{x_0, t_0 + (t + 1)}(M_{-1}^j) + \epsilon
\end{equation}
Hence if $t \in (-1 + \overline{r}/3, -1 + 2 \overline{r}/3)$, $t_0 \leq \overline{r}/3$, and $j > j_1$, then the above yield
\begin{equation}
F_{x_0, t_0}(M_t^j) < 1 + 3\epsilon
\end{equation}
Since $\epsilon$ was arbitrary, and $N$ has bounded second fundamental form, this is precisely what is need to apply White's Brakke regularity theorem for forced flows, theorem 4.1 in \cite{W}, to get uniform $C^{2, \alpha}$ bounds on $M_t^j$ for all $t \in (-1 + 4\overline{r}/9, -1 + 5\overline{r}/9)$ for some $\alpha \in (0,1)$. We can slightly change the $s_j$'s so that we instead have uniform $C^{2, \alpha}$ bounds on $M_{-1}^j$. 
Finally, observe that if $\Sigma_j$ is a sequence of closed surfaces converging to a closed surface $\Sigma_0$ as Radon measures and both the $\Sigma_j$'s and $\Sigma_0$ satisfy uniform $C^{2, \alpha}$ bounds, then the $\Sigma_j$'s must converge to $\Sigma_0$ uniformly in $C^2$. 
\end{pf}

Now, as in \cite{CM} we will construct a piece-wise MCF with a finite number of discontinuities that eventually becomes extinct in a round point, or one of the singularities encountered is noncompact (which could happen if, after doing some entropy decreasing perturbation, (1.1) fails to hold). This is by doing a smooth jump just before a (non-round) singular time, where we replace a time slice of the flow by a graph (in the normal bundle) over it. Moreover the perturbation we will show can be done so that the entropy decreases by at least a fixed $\epsilon' = \epsilon/4 > 0$ after each replacement, the $\epsilon$ of course coming from corollary 4.6 above. We repeat this until we get to a singular point where every tangent flow consists of shrinking spheres.  Recall again that after using the embedding theorems $M$ is a $K$-almost Brakke flow in $\R^\ell$. 
$\medskip$

By the assumption, all of the tangent flows at $t_0^{sing}$ are smooth, have multiplicity one (as described in the introduction), and correspond to compact self-shrinkers with diameter at most $D$ for some $D > 0$. In particular, since the tangent flows are compact, there is only one singular point $x_0 \in \R^\ell$. Let $\Sigma_0$ be a self-shrinker equal to the $t = -1$ time-slice of a multiplicity one tangent flow at the singularity. By the assumptions the lemma above gives a sequence $s_j > 0$ with $s_j \to 0$ so that 
\begin{equation}
\frac{1}{\sqrt{s_j}} M_{-s_j} \text{ converges in } C^2 \text{ to } \Sigma_0
\end{equation} 
There are two possibilities. First, if $\Sigma_0$ is the round sphere for every tangent flow at $x_0$, then the above implies that $M^0_{t_0^{sing} - s_j}$ is converging to a round sphere for every sequence $s_j \to 0$. 
Suppose instead that there is at least one tangent flow so that $\Sigma_0$ is not the round sphere. To proceed we will want to use our almost-monotonicity results so note since in (5.10) the sequence $s_j \to 0$, there is some $\overline{j}$ so that $K/\sqrt{s_{\overline{j}}}$ satisfies the assumptions of propositions 4.5 above for our choice of $\epsilon_0, \tau$. Let's determine what these should be before moving on. 
$\medskip$

After the rescaling (1.1) gives us that the diameter of $ \frac{1}{\sqrt{s_{\overline{j}}}} M_{t_0^{sing}-s_{\overline{j}}}$ is bounded by some $D < \infty$ as indicated above. Note that after rescaling $N$ by $s_j$ sufficiently small its curvatures (and derivatives thereof) as a submanifold of $\R^\ell$ can be made small enough so the sphere of radius $D$ will shrink to a point by using theorem 1.1 of \cite{H1}.  So after possibly taking $\overline{j}$ even larger by the typical maximum principle argument using the sphere of radius $D$ as a barrier then the time of the existence of its flow is bounded by say $\kappa$. So we take $\tau = \kappa$. 
$\medskip$

Since $\Sigma_0$ is not the round sphere by corollary 4.6 above (equation (1.1) let's us say the diameter of $\Sigma_0$ is bounded) we get a graph $\Gamma_0$ over $\Sigma_0$ with $\lambda(\Gamma_0) < \lambda(\Sigma_0) - \epsilon$ where $\epsilon > 0$ is a fixed constant given by the corollary. We set $\epsilon_0 = \epsilon/4$.
$\medskip$

With this in mind we relabel $M =  \frac{1}{\sqrt{s_{\overline{j}}}} M_{t_0^{sing}-s_{\overline{j}}}$. Note further rescalings by $S$, as long as $S > 1$, preserve $K < \delta$ needed for the choice of $\epsilon_0$ and $\tau$ above. Of course for large enough $j$, $\frac{1}{\sqrt{s_j}} > 1$. When $j$ is sufficiently large, $(\sqrt{s_j} \Gamma_0) + x_0$ is a graph in the normal bundle over $M_{t_0^{sing} - s_j}$ and 
\begin{equation} 
\lambda((\sqrt{s_j} \Gamma_0) + x_0) = \lambda(\Gamma_0) < \lambda(\Sigma_0) - \epsilon \leq \lambda( M_{t_0^{sing}-s_j} ) - \frac{3\epsilon}{4}
\end{equation}
where the first equality used the scale invariance of entropy and the last inequality used the almost-monotonicity of entropy under MCF. There's a small problem though, in that we see $(\sqrt{s_j} \Gamma_0) + x_0$ lies in $T_x N$, as a 3-plane in $\R^\ell$, but not necessarily $N$. So we want to project $(\sqrt{s_j} \Gamma_0) + x_0$ down to $N$. For this of course we need that entropy depends continuously on $\Sigma$ in at least certain cases so our work in section 5 above comes in handy. 
$\medskip$

Recalling proposition 5.4 above, to proceed then we need to show for large enough rescalings of $\R^\ell$ the projection of $(\sqrt{s_j} \Gamma_0) + x_0$ back onto $N$, which we'll denote $\widetilde{\Gamma}$, is as close as we want to $(\sqrt{s_j} \Gamma_0) + x_0$ in $C^{2,\alpha}$ topology. More precisely:
\begin{lem} For each $D, \rho > 0$, there exists $\xi > >0$ so that $\xi(N) \cap B_{D}(x)$ is a graph of a function $f$ over $T_xN$ and $||f||_{C^{3}}< \rho$.
\end{lem}
\begin{pf} 
The second fundamental form $A_N$ of $N$ is bounded in norm initially and scales by $A \to \frac{1}{\xi^2} A$ under rescaling by $\xi$ so taking $\xi$ as large as we want we make $|A_N|^2$ as small as we want, which implies $N$ is graphical - see lemma 2.4 in \cite{CM2}. For a graph, $C^{2,\alpha}$ bounds scale by inverse power so that, possibly taking $\xi$ larger, we can arrange so that $||f||_{C^3}< \rho$ for a given $\rho > 0$. 
\end{pf} 
$\medskip$

We see then, possibly taking $\overline{j}$ above larger, without loss of generality $|\lambda((\sqrt{s_j} \Gamma_0) + x_0) - \lambda(\widetilde{\Gamma})| < \epsilon/4$. Hence:
\begin{equation} 
\lambda(\widetilde{\Gamma}) \leq \lambda( M_{t_0^{sing}-s_j} ) - \frac{\epsilon}{2}
\end{equation}
We then replace $ M_{t_0^{sing}-s_j} $ with $\widetilde{\Gamma}$ and restart the flow (note its diameter will be less than $D$ so the time of existence of the flow less than $\kappa$ above).
$\medskip$

Note since we already rescaled space by $\frac{1}{\sqrt{s_{\overline{j}}}}$ above, before the first replacement, so that the almost monotonicity lemmas and ``flatness of $N$" along the lines of lemma 6.3 hold with $\epsilon_0 = \epsilon/4$, $\tau = \kappa$ (again, as long as later rescalings by say $S$ are so that $S > 1$ - this is done by throwing out $s_j$ (in subsequent blowup sequences) with $s_j > 1$). The entropy by the next singularity can only go up by an additional $\epsilon/4$, which is to say the entropy of the tangent flow $\Sigma_1$ at the ``next" singular time satisfies:
\begin{equation}
\lambda(\Sigma_1) < \lambda(\Sigma_0) - \frac{\epsilon}{4}
\end{equation}
Since the entropy of $M$ was initially finite and entropy goes down a uniform constant $\epsilon/2 > 0$ at each replacement and is guaranteed not to increase by more than $\epsilon/4$ under the flow to the next time, this can only occur a finite number of times until all tangent flows are round.
$\medskip$

\section{Concluding remarks.}
One sees that generally speaking given a hypersurface $M$ so that equation (1.1) holds up to the first singular time (1.1) can't be expected to hold after the perturbation above; it is part of the theorem (that (1.1) holds for the whole piecewise flow) that this is so for us. For example, if a hypersurface $M \subset N$ shrinks to a point and is not genus 0, we see that (1.1) must eventually fail to hold for some perturbation because the genus of the surface is unchanged after each replacement and preserved under the flow. Indeed, a corollary of theorem 1.1 is that a piecewise MCF starting out of a positive genus surface will eventually encounter a non compact (after rescaling) singularity. 
$\medskip$

 It follows by Brendle's classification in \cite{Bre}, shown after \cite{CM} was completed, of genus 0 closed embedded self-shrinkers in $\R^3$ as precisely the round sphere of radius 2 that conversely if $M$ is genus zero and shrinks to point that it must do so to a round point so in that sense our theorem is weaker than what is already known to be true. On the other hand, our methods may readily be generalized to higher dimensions save for the fact that, as pointed out above, we don't have quite as good compactness theorems for self shrinkers in higher dimensions so the ``off the shelf" possible theorems aren't quite as strong -one might possibly for example impose a type 1 blowup rate assumption along with equation 1.1 and use that self shrinkers with bounded curvature and diameter are a compact set. We also feel it is interesting to understand how to apply some of Colding and Minicozzi's theory to the curved setting. 
$\medskip$

On that note in $\R^{3}$, there is at least one known example of a positive genus closed self shrinker (which of course shrink to a point), the Angenant torus \cite{Ang} (see \cite{P} for higher dimensional analogues). Also relevant (although they don't shrink to points) are the more recent noncompact self shrinkers by Kapouleas, Kleene, and Moller in \cite{KKM} of higher genus obtained by gluing or those obtained by Ketover in \cite{K} using min-max methods (in fact, the Ketover examples are possibly actually compact but that hasn't been decided). To the author's knowledge though, it is unknown if there are examples in general ambient manifolds of hypersurfaces $\Sigma$ of positive genus that shrink to points under the flow (convex enough surfaces, which are spheres, will by \cite{H1}); one imagines it could be possible to get examples using an inverse function theorem argument by ``importing" examples like Angenant torus. Related to this, one wonders if there is a good ``intrinsic" notion of self shrinker on an ambient manifold $N$ that would help facilitate this. 
$\medskip$

It is also worth pointing out that singularities forming as $t \to \infty$ $\textit{might}$ possibly occur under the flow in the sense that the surface could flow to a minimal cone as $t \to \infty$, so that's why we stipulate above explicitly in the statement of theorem 1.1 that $t_0 < \infty$. For example, there are examples by Valasquez, see \cite{V}, of surfaces that develop the Simons cone as a singularity.


\begin{thebibliography}{9}

\bibitem{Ang}  Angenent, Sigurd. Shrinking doughnuts.Progress in Nonlinear Differential Equations and their Applications, 7, 21-38.

\bibitem{And} Andrews, Ben; Li, Haizhong; Wei, Yong. $\mathcal{F}$-stability for self-shrinking solutions to mean curvature flow. Asian J. Math. 18 (2014), no. 5, 757--778. 

\bibitem{BW} Bernstein, Jacob; Wang, Lu. A Sharp Lower Bound for the Entropy of Closed Hypersurfaces up to Dimension Six. arXiv:1406.2966

\bibitem{B} Brakke, Kenneth. The Motion of a Surface by its Mean Curvature. Princeton University Press, 1978.

\bibitem{Bre} Brendle, Simon. Embedded self-shrinkers of genus 0. arXiv:1411.4640 

\bibitem{CM} Colding, Tobias; Minicozzi II, William. Generic mean curvature flow I; generic singularities. Annals of Mathematics. (2) 175 (2012), 755-833. 

\bibitem{CM1} Colding, Tobias; Minicozzi II, William. Smooth compactness for self-shrinkers. arXiv:0907.2594

\bibitem{CM2} Colding, Tobias; Minicozzi II, William. A course in minimal surfaces. Graduate studies in mathematics; v. 121. 

\bibitem{H} Huisken, Gerhard. Asymptotic behavior for singularities of the mean curvature flow. J. Differential Geom. 31 (1990), no. 1, 285--299. 

\bibitem{H1} Huisken, Gerhard. Contracting convex hypersurfaces in Riemannian manifolds by their mean curvature. Inventiones mathematicae 84 (1986): 463-480.

\bibitem{I} Ilmanen, Tom. Singularities of the mean curvature flow of surfaces. Preprint.

\bibitem{KKM} Kapouleas, Nicos; Kleene, Stephen; Moller, Niels.  Mean curvature self-shrinkers of high genus: Non-compact examples. J. Reine Angew. Math. (2012), to appear.

\bibitem{K} Ketover, Daniel. Self-shrinking platonic solids. arXiv:1602.07271

\bibitem{Mant} Mantegazza, Carlos. Lecture Notes on Mean Curvature Flow. Springer Basel, 2013. 

\bibitem{N} Nash, John. The imbedding problem for Riemannian manifolds. Annals of Mathematics. (1) 63 (1956), 20-63.

\bibitem{P} McGrath, Peter. Closed Mean Curvature Self-Shrinking Surfaces of Generalized Rotational Type. arXiv:1507.00681

\bibitem{W} White, Brian. Stratification of minimal surfaces, mean curvature flows, and harmonic maps. Journal für die reine und angewandte Mathematik 488 (1997): 1-36. 

\bibitem{W1} White, Brian. A local regularity theorem for mean curvature flow.  Annals of Mathematics. (2) 161 (2005),  1487-1519. 

\bibitem{V} Velázquez, J. J. L.. Curvature blow-up in perturbations of minimal cones evolving by mean curvature flow. Annali della Scuola Normale Superiore di Pisa - Classe di Scienze 21.4 (1994): 595-628.

\bibitem{YL} Lee, Yng-Ing. Lue, Yang-Kai. The stability of self-shrinkers of mean curvature flow in higher codimension. Trans. Amer. Math. Soc. 367 (2015), 2411-2435.

\bibitem{Z} Zhu, Jonathan. On the entropy of closed hypersurfaces and singular self-shrinkers. arXiv:1607.07760
\end{thebibliography}
\end{document}